\documentclass{amsart}

\usepackage{amssymb}

\newtheorem{thm}{Theorem}

\newtheorem{conj}[thm]{Conjecture}
   
\theoremstyle{definition}
\newtheorem{defn}[thm]{Definition}

\newtheorem{say}[thm]{}
\newtheorem{exmp}[thm]{Example}

\newtheorem{exrc}[thm]{Exercise}

\newtheorem{ques}[thm]{Question}    

\newtheorem{aside}[thm]{Aside}          
\newtheorem{ack}{Acknowledgments}

\newtheorem{defn-thm}[thm]{Definition--Theorem}  
\newtheorem{defn-lem}[thm]{Definition--Lemma}  

\theoremstyle{remark}

\renewcommand{\c}[0]{{\mathbb C}}  

\renewcommand{\o}[0]{{\mathcal O}} 
\newcommand{\z}[0]{{\mathbb Z}}

\newcommand{\p}[0]{{\mathbb P}}
\newcommand{\f}[0]{{\mathbb F}}
\newcommand{\q}[0]{{\mathbb Q}}
\newcommand{\map}[0]{\dasharrow}
\newcommand{\qtq}[1]{\quad\mbox{#1}\quad}

\newcommand{\pic}[0]{\operatorname{Pic}}

\newcommand{\codim}[0]{\operatorname{codim}}

\newcommand{\aut}[0]{\operatorname{Aut}}

\newcommand{\jac}[0]{\operatorname{Jac}}

\newcommand{\univ}[0]{\operatorname{univ}}

\newcommand{\mor}[0]{\operatorname{Mor}} 
 
\newcommand{\tsum}[0]{\textstyle{\sum}}

\newcommand{\curves}[0]{\operatorname{Curves}}

\def\into{\DOTSB\lhook\joinrel\to}

\begin{document}
\bibliographystyle{amsalpha}

\title[Curves on varieties]{Holomorphic and pseudo-holomorphic curves\\ on 
rationally connected varieties}
\author{J\'anos Koll\'ar}
\begin{abstract} These notes give a short introduction
to the study of curves on algebraic varieties.
The main emphasis is on families of genus $0$  curves.
After an elementary proof of the dimension formula
for the space of curves, we summarize the basic
properties of uniruled and of  rationally connected varieties.
The last section is devoted to a conjectural
characterization of rationally connected varieties using  symplectic geometry.
\end{abstract}


\maketitle

One of the aims of algebraic geometry is to understand all subvarieties of a 
given algebraic variety. For simplicity, assume that $X$ is a smooth
projective variety over $\c$. We have a very successful -- by now classical --
theory that describes codimension 1 subvarieties, that is,
hypersurfaces or divisors on $X$.

\begin{say}[Codimension 1 theory]\label{codim1.say}
(For details, see, for instance \cite[Chap.1]{gri-har}.)

 Each codimension 1 subvariety
$Z\subset X$ determines a line bundle (equivalently,
 a rank 1 locally free sheaf)
 $L:=\o_X(Z)$ and a section
$s\in H^0(X,L)$; the constant 1 section of $\o_X(Z)$.
($L$ is determined uniquely, $s$ is determined  
up to a multiplicative constant.)
It is thus sufficient to describe the pairs $(L,s)$.
We proceed in 3 steps.
\medskip

(\ref{codim1.say}.1) As a topological line bundle, $L$ is determined
by its Chern class $c_1(L)\in H^2(X(\c),\z)$.
Furthermore, a cohomology class $\alpha\in H^2(X(\c),\z)$ is the Chern class
of  a holomorphic/algebraic line bundle 
iff the image of $\alpha$ in $H^2(X(\c),\c)$ lies in
$H^{1,1}(X(\c),\c)$. 
\medskip

(\ref{codim1.say}.2) The set of all 
holomorphic/algebraic line bundles with a given 
 $c_1(L)\in H^2(X(\c),\z)$ is parametrized by
(more precisely, is a principal homogeneous space under)
an Abelian variety, called the Picard variety, $\pic^0(X)$.
It is biholomorphic to 
$H^{0,1}\bigl(X(\c),\c\bigr)/H^1\bigl(X(\c),\z\bigr)$. 

\medskip

(\ref{codim1.say}.3) The sections of a line bundle $L$ form a  vector space 
$H^0(X,L)$, and nonzero sections up to a multiplicative constant
are parametrized by the corresponding projective space
$|L|:=\p\bigl(H^0(X,L)^{\vee}\bigr)$.
In  terms of \v{C}ech cohomology, computing $H^0(X,L)$ is essentially
a linear algebra problem. On the other hand,
the  Hirzebruch-Riemann-Roch theorem and Serre's vanishing theorem say
 that, for sufficiently ample
$L$, the dimension of $H^0(X,L)$ can be computed
from the Chern classes, thus from topological data.
\end{say}

{\bf Summary:} The topology of $X(\c)$ and Hodge theory 
determine Steps (\ref{codim1.say}.1--2)
and linear algebra governs Step (\ref{codim1.say}.3).
\medskip

Ideally, we would like to have a similar description of
higher codimension subvarieties, but this goal is very far off.

The aim of these notes is to focus on the next simplest case,
the study of 
1-dimensional subvarieties, that is, curves in $X$.
Despite Poincar\'e duality between divisors and curves,
 the study of curves turns out to be quite a bit harder
than the theory of divisors. There are
 still many  deep open questions.

Let us start with the analog of (\ref{codim1.say}.1).

\section{Homology classes of curves}

Given any curve $C\subset X$, it has a homology class
$[C]\in H_2(X(\c),\z)$. Our first question is: Which
homology classes can be realized by algebraic curves?
As before, this class has Hodge type $(1,1)$. 
We usually think of Hodge theory as living on
cohomology groups, so we identify 
$H_2(X(\c),\c)$ with $H^{2n-2}(X(\c),\c)$ (where $n=\dim X$),
and think of
$[C]$ as an integral cohomology  class of type $(n-1,n-1)$.

The hard Lefschetz theorem \cite[p.122]{gri-har}
implies that if $L$ is an ample divisor then
$$
H^{n-1,n-1}\bigl(X(\c),\c\bigr)=L^{n-2}\cdot H^{1,1}\bigl(X(\c),\c\bigr).
$$
In particular, a rational multiple of every curve class
can be written as $\bigl(L^{n-2}\cdot D\bigr)$ for some 
(not necessarily effective) divisor $D$ on $X$.
This  gives the complete answer with rational
coefficients. The problem is, however, more subtle over $\z$.

To illustrate how little is known,
let us discuss two conjectures about
homology classes of curves.

Consider smooth hypersurfaces
$X^n_e\subset \p^{n+1}$ of degree $e$.
By another Lefschetz theorem \cite[p.156]{gri-har}, the natural map
$H_2(X(\c),\z)\to H_2(\c\p^{n+1},\z)$ is an isomorphism
if $n\geq 3$. Thus there is a class $\ell\in H_2(X(\c),\z)$
whose image in $H_2(\c\p^{n+1},\z)$ is the homology class of a  line
and $H_2(X(\c),\z)=\z[\ell]$.

The original form of the Hodge conjecture implies that
$\ell$ is the homology class of a ($\z$-linear combination of)
algebraic curves, while the Noether-Lefschetz theorem suggests
that all curves on $X$ may be complete intersections with
surfaces.

\begin{conj} \cite{g-h-85} Let $X^n\subset \p^{n+1}$ be a  very general
smooth hypersurface and $C\subset X$ an algebraic curve.
Assume that $\deg X$ is large enough
(maybe $\deg X\geq 2n$  is sufficient).
Then $[C]\in H_2(X(\c),\z)=\z[\ell]$
is a multiple of $\deg X\cdot [\ell]$.
\end{conj}

A straightforward dimension count shows that every
hypersurface $X^n\subset \p^{n+1}$ 
contains a line if $\deg X\leq 2n-1$,
but the general hypersurface does not contain a line if
$\deg X\geq 2n$. (This, however, does not exclude the possibility
that it contains a $\z$-linear combination of
algebraic curves whose degree is 1.)

An easy result \cite{k-trento} shows that if
$(6,k)=1$, 
 $X$ is a very general hypersurface of degree $3k^2$
and $C\subset X$ is any curve, then $[C]$
is a multiple of $k\cdot [\ell]$.

By contrast,  for every degree there are smooth hypersurfaces
that contain a line. Thus topology and Hodge theory together 
do not give enough information to decide which 2-dimensional homology classes
are algebraic.

The opposite may hold for  smooth projective Fano varieties
(that is, varieties with $-K_X$  ample) or more generally
for rationally connected varieties (\ref{RCV.defn}). 
In the Fano case, these conjectures go back to Fano and Iskovskikh;
the general case is proposed in \cite{voi-1}.

\begin{ques}\label{line-fano.conj}
Let $X$ be a smooth projective  variety.
\begin{enumerate}
\item  Is $ H_2(X(\c),\z)$  generated by the homology classes of 
rational curves if $X$ is Fano or, more generally, if $X$ is 
 rationally connected?
\item Assume that  $X$ is Fano  and 
$\eta\in H_2(X(\c),\z)$ 
such that $m\eta$ is the  homology class of an effective algebraic curve
for some $m>0$. 
  Is $\eta$  the homology class of an effective algebraic curve, 
all of whose irreducible 
components are rational?
\end{enumerate}
\end{ques}

Both of these are known with $\q$-coefficients.
(\ref{line-fano.conj}.1) is quite easy  \cite[IV.3.13]{rc-book}
and (\ref{line-fano.conj}.2) follows from Mori's cone theorem
\cite{mori-cone}. 

In dimension 3  the positive answer to 
(\ref{line-fano.conj}.1) is   a special case
of  \cite{voi-2}. The original conjecture of Fano and Iskovskikh
asks (\ref{line-fano.conj}.2) in case the second Betti number of $X(\c)$ is 1.
In dimension 3 this is established as part of the classification
of Fano 3-folds.

If $X$ is rationally connected, the technique of ``deformation of combs''
(cf.\ \cite[Sec.II.7]{rc-book})  shows that
the subgroup of $H_2(X(\c),\z)$ generated by
homology classes of curves (resp.\ rational curves)
is invariant under smooth deformations of $X$.

\begin{exmp} 1. Let $\pi:X\to \p^n$ be a smooth degree $d$ cover of $\p^n$,
given by an affine equation of the form $z^d=f_{dr}(x_1,\dots,x_n)$
where $f_{dr}$ had degree $dr$. The branch locus $Z\subset \p^n$ is 
(the projective closure of) the hypersurface $(f_{dr}=0)$.

For special choices of $f_{dr}$, there are lines
$L\subset \p^n$ that are $d$-fold tangent to $Z$ at $r$-points
$P_1,\dots, P_r$.
Then $\pi^{-1}(L)$ splits as a union of $d$ irreducible curves,
each mapping isomorphically to $L$.
Let $L_i\subset X$ be any of these curves and $P_j\in X$
any of the points. 
Set $Y:=B_{P_j}X$, the blow-up of $X$ at the point $P_j$
with exceptional divisor $E$.
Let $\alpha\in H_2(Y(\c),\z)$ denote the 
homology class of the birational transform $L'_i\subset Y$ of $L_i$.
The class $\alpha$ is characterized by the properties
$\alpha\cdot E=1$ and $\alpha\cdot \pi_Y^*H=1$
where $H$ is the hyperplane class on $\p^n$ and
$\pi_Y:Y\to \p^n$ the composite of the blow-down map with $\pi$. 
Note that the class $\alpha$ is well defined for any
$B_PX$, no matter where the point $P\in X$ lies.

An easy dimension count as in \cite[IV.2.12]{rc-book} shows that
if $(d-1)r\geq n$ then there is no such line for
general $P\in Z$. Moreover, if $r\geq 2$, then we can choose
$P$ and $Z$ such that there is a line $L_P$ that is
 $d$-fold tangent to $Z$ at $P$, but every such line 
has transverse intersection with $Z$ 
at some other point.

Then $\pi^{-1}(L_P)$ is irreducible and has multiplicity $d$ at $P$.
Hence its birational transform on $B_PX$ represents the
homology class $d\alpha$. However, there are no
effective curves in the homology classes
$\alpha,\dots, (d-1)\alpha$.

If $(d-1)r= n$ then $X$, and hence $B_PX$,  are rationally connected.
$X$ is Fano but $B_PX$ is not Fano.

2. By \cite{art-mum}, there is a smooth, projective,
unirational 3-fold $X$ such that $H_2(X(\c), \z)$ contains a
2-torsion element $\alpha$. 
One can write $\alpha=[C_1]-[C_2]$ where the $C_i\subset X$
are smooth rational curves.

Using \cite{dJ-brauer}, one can show that 
there is a smooth 4-fold $\pi:Y\to X$ that is analytically locally 
a $\p^1$-bundle but it has no rational section.
Thus every line bundle on $Y$ is of the form
$\pi^*L(mK_Y)$ where $L$ is a line bundle on $X$.

Let $C\cong \p^1$ be a fiber of $\pi$, then
$\bigl(C\cdot c_1(\pi^*L(mK_Y))\bigr)=-2m$, hence even.
Since $h^2(Y,\o_Y)=0$, we see that every class in $H^2(X(\c), \z)$
is algebraic, thus $(C\cdot L)\in 2\z$ for every 
$L \in H^2(X(\c), \z)$. By Poincar\'e duality, there is a
homology class $\beta\in H_2(X(\c), \z)$ such that $[C]=2\beta$
but $\beta$ is not the homology class of an effective algebraic curve.
Nonetheless, one can find liftings $C'_i\subset Y$ of the $C_i$
such that $\beta=[C'_1]-[C'_2]$. Thus
(\ref{line-fano.conj}.2) does not hold for $Y$ but
(\ref{line-fano.conj}.1) does.

In this example $Y$ is unirational but  not Fano.
However, I see no reason why similar Fano examples
should not exist. Thus the answer to (\ref{line-fano.conj}.2)
may well be negative.
\end{exmp}


\begin{exrc} \cite{ful-etal} Let $X$ be a proper toric variety.
Then $H_2(X(\c),\z)$ is generated by the
homology classes of 1-dimensional orbits.
\end{exrc}

\section{The dimension of the space of curves}{\ }

Here we consider the analogs of
(\ref{codim1.say}.2--3) together for curves,
since it is not known how to separate the Abelian part
(\ref{codim1.say}.2) from the linear algebra part
(\ref{codim1.say}.3). Actually, it seems quite unlikely that
considering the two parts separately would make sense for
families of curves on varieties of dimension $\geq 3$.
 Our final results are  (\ref{dim.est.gen.thm})
and (\ref{dim.chow/hilb.thm}); we build up to them 
through a series of examples.

\begin{exmp}[Maps of rational curves to $\p^n$] \label{dc.exmp.1}
Giving a morphism $f:\p^1\to \p^n$ whose image
has degree $d$ is equivalent to giving
$n+1$ homogeneous polynomials $(f_0,\dots, f_n)$ of degree $d$
 without common root (up to a multiplicative constant).
Thus the space of all such maps is an open subset
of a projective space of dimension $(n+1)(d+1)-1$:
$$
\mor_d(\p^1,\p^n)\subset \p^{(n+1)(d+1)-1}.
$$
\end{exmp}

\begin{exmp}[Maps of rational curves to hypersurfaces]
 \label{dc.exmp.2} Let $X^n_e\subset \p^{n+1}$ be a hypersurface
of degree $e$ given by an equation
$G(x_0,\dots, x_{n+1})=0$.  
The image of a morphism $f:\p^1\to \p^{n+1}$ as in (\ref{dc.exmp.1}) lies in
$X$ iff $G(f_0,\dots, f_{n+1})\equiv 0$. The latter is a
homogeneous polynomial of degree $de$, so its vanishing is
equivalent to $de+1$ equations in the coefficients of the $f_i$.
Thus the space of all degree $d$ maps $\p^1\to X$
is either empty or a closed subset
of  $\mor_d(\p^1,\p^{n+1})$ of codimension $\leq de+1$.
In particular
$$
\dim_{[f]} \mor_d(\p^1,X^n_e)\geq (n+2)(d+1)-1-(de+1)=(n+2-e)d+n.
$$
(There is no reason to assume that  $\mor_d(\p^1,X^n_e)$
is pure dimensional, and $\dim_{[f]}$ denotes the dimension
at the point $[f]\in \mor_d(\p^1,X^n_e)$ corresponding to $f$.
Actually we proved something stronger: the lower bound holds
for every irreducible component
containing $[f]$.)

With hindsight masquerading as prescience, we can write the
right hand side as
$$
\dim_{[f]} \mor_d(\p^1,X^n_e)\geq c_1(X)\cdot f_*[\p^1] +\dim X.
\eqno{(\ref{dc.exmp.2}.1)}
$$
\end{exmp}

\begin{exmp}[Maps of  curves to $\p^n$] \label{dc.exmp.3} 
Now let us replace $\p^1$ by any curve $C$.
To start, one can think of $C$ as a smooth projective curve, but all
the arguments  work if $C$ is  a 1-dimensional, smoothable 
 projective scheme over a field $k$ such that
$H^0(C,\o_C)\cong k$.

First, let us fix a line bundle $L$ on $C$.
Giving a morphism $f:C\to \p^n$ 
such that $f^*\o_{\p^n}(1)\cong L$
 is equivalent to giving
$n+1$ section 
$f_0,\dots, f_n\in H^0(C,L)$ 
without common root (up to a multiplicative constant).
Thus the space of all such maps is an open subset
of a projective space of dimension $(n+1)h^0(C,L)-1$,
hence
$$
\dim \mor_L(C,\p^n)= (n+1)h^0(C,L)-1,
$$
where $\mor_L$ denotes those morphisms for which 
$f^*\o_{\p^n}(1)\cong L$.
\end{exmp}

\begin{exmp}[Maps of curves to hypersurfaces] \label{dc.exmp.4} 
 Now let $X^n_e\subset \p^{n+1}$ be a hypersurface
of degree $e$ given by an equation
$G(x_0,\dots, x_{n+1})=0$.  
Then $G(f_0,\dots, f_{n+1})$ is a section of 
$L^e$, hence its identical vanishing imposes
$ h^0(C,L^e)$ conditions. 
As in (\ref{dc.exmp.2}),  the space of all such maps $C\to X$
is either empty or a closed subset
of  $\mor_L(C,\p^{n+1})$ of codimension $\leq h^0(C,L^e)$.
In particular
$$
\dim_{[f]} \mor_L\bigl(C,X^n_e\bigr)\geq (n+2)h^0(C,L)-1- h^0(C,L^e).
\eqno{(\ref{dc.exmp.4}.1)}
$$
If $\deg L$ is large enough, then $h^0(C,L)=\deg L+1-g$
and we can rewrite the formula as
$$
\dim_{[f]} \mor_L\bigl(C,X^n_e\bigr)\geq c_1\bigl(X^n_e\bigr)\cdot f_*[C] +
n\cdot \chi(\o_C)-g(C).
\eqno{(\ref{dc.exmp.4}.2)}
$$
This fits very nicely with (\ref{dc.exmp.2}.1),
 except for the $-g(C)$ term at the end.
Remember now that at the beginning we fixed not just $\deg L$ but $L$ itself.
All line bundles of  degree $d$ are parametrized by
the $g(C)$-dimensional $\jac_d(C)$, which is   (non-canonically) isomorphic to
the Jacobian of $C$
(see, for instance, \cite[Sec.2.7]{gri-har}).
 Thus, letting $L$ vary 
should increase the dimension by
$g(C)$ and we should get that 
$$
\dim_{[f]} \mor(C,X)\geq c_1(X)\cdot \phi_*[C] +\dim X\cdot \chi(\o_C).
\eqno{(\ref{dc.exmp.4}.3)}
$$
(Note, however, that  our results included a caveat that
the spaces may be empty. Thus, as we vary $L$
in $\jac_d(C)$, we could get empty spaces for some $L$ where we expect
something  positive dimensional.  Hence (\ref{dc.exmp.4}.3)
is not yet firmly proven.)
\end{exmp}

\begin{exmp}[Curves on $\p^n$]\label{philos.say.1}
So far we have been estimating the dimension of the space
of morphisms from a fixed curve to $\p^n$
with fixed $f^*\o_{\p^n}(1)$. In order to estimate
the dimension of the space of curves on $\p^n$, we need to take three
differences into account. 
(If $h^1(C,\o_C)=0$ then only the second and third steps are needed.)

First, we need to work out precisely how to let
$f^*\o_{\p^n}(1)$ vary in $\jac_d(C)$.

Second, if $g=0$ or $g=1$,
then every curve has many automorphisms. Thus a single rational curve 
$C\subset \p^n$ gives rise to a whole family of 
birational maps $\p^1\to C$ parametrized by  $\aut(\p^1)=PGL(2)$.
Similarly, every elliptic curve 
$C\subset \p^n$ gives rise to a 1-dimensional family of maps.

Third, for $g\geq 1$, we can vary the curve $C$ in its moduli space.

To do the first, let $\jac_d(C)$ denote the family of degree
$d$ line bundles on $C$.
If $d\geq 2g-1$, the universal $H^0$ gives a vector bundle
${\mathcal H}_d\to \jac_d(C)$.
As before, the degree $d$ maps from $C$ to
$\p^n$ are parametrized by an open subset of the
projective space bundle
$$
\p_{\jac_d(C)}\Bigl( \bigl({\mathcal H}_d^{\vee}\bigr)^{n+1}\Bigr)
\to \jac_d(C).
$$
This gives the expected dimension formula
$$
\begin{array}{rcl}
\dim_{[f]} \mor_d(C,\p^n)&\geq & c_1(\p^n)\cdot f_*[C] +
n\cdot \bigl(1-g(C)\bigr)\\
 & = & (n+1)\deg C + n\cdot \bigl(1-g(C)\bigr).
\end{array}
\eqno{(\ref{philos.say.1}.1)}
$$
For $g\geq 2$, consider a $(3g-3)$-dimensional family  ${\mathcal M}_g$ 
of curves
of genus $g$ (in practice, some finite cover of the moduli space)
over which there is a universal family of curves
${\mathcal C}_g\to {\mathcal M}_g$.
(If $C$ is singular, we use the local deformation space
of $C$ as in \cite[II.1.11]{rc-book}.)
Let ${\mathcal J}_d\to {\mathcal M}_g$ denote the universal family
of degree $d$ components of the Jacobians.
As before, if $d\geq 2g-1$, the universal $H^0$ gives a vector bundle
$$
{\mathcal H}_d\to{\mathcal J}_d\to  {\mathcal M}_g,
$$
and $\dim {\mathcal H}_d=(3g-3)+g+(d+1-g)$.
Thus the degree $d$ maps from  genus $g$ curves to
$\p^n$ are parametrized by an open subset of the
projective space bundle
$$
\p_{{\mathcal J}_d}\Bigl( \bigl({\mathcal H}_d^{\vee}\bigr)^{n+1}\Bigr)
\to{\mathcal J}_d\to  {\mathcal M}_g.
$$
Hence, if $C\subset \p^n$ is a smooth curve 
of degree $\geq 2g(C)-1$, then
$$
\begin{array}{rcl}
\dim_{C} \curves(\p^n)&\geq & c_1(\p^n)\cdot [C] +
(n-3)\cdot \bigl(1-g(C)\bigr)\\
 & = & (n+1)\deg C + (n-3)\cdot \bigl(1-g(C)\bigr),
\end{array}
\eqno{(\ref{philos.say.1}.2)}
$$
where $\curves(\p^n)$ denotes either the
Chow variety or the Hilbert scheme of curves on $\p^n$.
(These two spaces agree near smooth or normal subvarieties,
\cite[Cor.I.6.6.1]{rc-book}.)
It is easy to check that the formulas also work if $g(C)\leq 1$.
\end{exmp}

\begin{exmp}[Maps of curves to  varieties]
Let $X^n\subset \p^N$ be a smooth projective subvariety.
In general $X$  needs many defining  equations
$G_i=0$. Correspondingly, 
the image of a morphism $f:C\to \p^N$  lies in
$X$ iff $G_i(f_0,\dots, f_N)\equiv 0$ for every $i$.
This implies that $\mor_L(C,X)$ is a closed algebraic subvariety
of $\mor_L(C,\p^N)$, but we can not estimate its codimension
unless we know the degrees of  the equations $G_i$.

If $X^n\subset \p^N$ is a complete intersection, 
there is a ``natural'' choice for the $G_i$ and
everything we said before 
works out.
In particular, (\ref{dc.exmp.4}.2--3) still hold.
In general, however, we get an estimate that is much worse
and depends on the choice of the equations $G_i$.

One can, however, easily compute from this  presentation
the tangent spaces of $\mor_L(C,X)$.
\end{exmp}

\begin{thm} \label{tangent.space.thm}
Let $X$ be a smooth   variety, $C$  a
proper  curve 
such that $H^1(C, \o_C)=0$
and $f:C\to X$ a morphism.
Then
 $$
T_{[f]}\mor(C,X)=H^0\bigl(C, f^*T_X\bigr).
$$
\end{thm}

Proof. Assume  that $X^n\subset \p^N$ is  given by the equations
$G_i=0$ and let
$f:C\to X$ be given by
$(f_0:\dots : f_N)$ where the $f_i$ are sections of
$f^*\o_X(1)$. Inside $\mor(C,\p^N)$ the
tangent directions are given by the deformations
$\bigl(f_0+th_0:\dots : f_N+th_N\bigr)$  where the $h_i$ are also sections of
$f^*\o_X(1)$.
The corresponding tangent vector is in the tangent space
of $\mor(C,X)$ iff
$$
G_i\bigl(f_0+th_0,\dots, f_N+th_N\bigr)\equiv 0 \mod (t^2)
\qquad \forall\ i.
$$
Using the Taylor expansion, this is equivalent to
$$
{\tsum_j} \frac{\partial G_i}{\partial x_j}\cdot h_j=0\qquad \forall\ i.
$$
The latter holds if 
$\bigl(h_0,\dots, h_N\bigr)$ maps to a section of $f^*T_X
\subset f^*T_{\p^N}$
in the exact sequence
$$
0\to f^*\o_{\p^N}\to f^*\o_{\p^N}(1)^{N+1}\to f^*T_{\p^N}\to 0.
$$
Since $f^*\o_{\p^N}=\o_C$, by taking cohomology we get
the exact sequence
$$
H^0\bigl(C, f^*\o_{\p^N}(1)\bigr)^{N+1}\to 
H^0\bigl(C,f^*T_{\p^N}\bigr)\to
H^1\bigl(C,\o_C\bigr)=0,
\eqno{(\ref{tangent.space.thm}.1)}
$$
thus every section of $H^0\bigl(C, f^*T_X\bigr)$ 
is the image of some $\bigl(h_0,\dots, h_N\bigr)$.
The same argument shows that if $g(C)\neq 0$ then
 $$
T_{[f]}\mor_L(C,X)=
\ker\Bigl[H^0\bigl(C, f^*T_X\bigr)\to H^1\bigl(C,\o_C\bigr)\Bigr],
\eqno{(\ref{tangent.space.thm}.2)}
$$
where the map on the right comes from the sequence
(\ref{tangent.space.thm}.1).\qed

\medskip
While we usually think of an algebraic variety as sitting inside
a larger dimensional projective space, one can also represent an
n-dimensional variety as a finite branched cover  $\pi:X^n\to \p^n$.
This was Riemann's original point of view   for algebraic curves.
In higher dimensions it can be obtained by repeatedly projecting
$X^n\subset \p^N$ from points outside it,
until we get a dominant morphism $\pi:X^n\to \p^n$.
This will allow us to prove the expected lower bounds for
the spaces of maps.

\begin{say}[Lifting deformations to branched covers]\label{lifting.say}
Let $\pi:X\to \p^n$ be a finite surjection 
with ramification divisor $R\subset X$
and $f:C\to X$ a morphism from a reduced curve to $X$.
Assume the following genericity conditions:
\medskip

{\it Assumptions} \ref{lifting.say}.1.
\begin{enumerate}
\item[i)]  $C$ is smooth and $\pi\circ f:C\to \p^N$ is unramified at
 every $p\in f^{-1}(R)$,
\item[ii)]  near $f(C)\cap R$, $R$ is smooth and the ramification index
of $\pi$ is 2, 
\item[iii)]   
$R$  and $f(C)$ 
 intersect transversally.
\end{enumerate}
\medskip

Equivalently, we can choose local analytic coordinates
$(x_1,\dots, x_n)$ near $p$ and $(y_1,\dots, y_n)$ near $\pi(p)$
such that $\pi$ is  given by
$$
\pi: (x_1,\dots, x_n)\mapsto  (x_1,\dots, x_{n-1}, x_n^2), 
$$
and $C$ is parametrized as
$$
f: t\mapsto \bigl(f_1(t),\dots, f_n(t)\bigr),
$$
where $f_i(0)=0$, $\bigl(f'_1(0),\dots, f'_{n-1}(0)\bigr)\neq (0,\dots,0)$
and $f'_n(0)\neq 0$. Thus $\pi\circ f$ is given by
$$
\pi\circ f: t\mapsto \bigl(f_1(t),\dots,  f_{n-1}(t), f^2_n(t)\bigr)
=: \bigl(g_1(t),\dots, g_n(t)\bigr).
\eqno{(\ref{lifting.say}.2)}
$$
Its image is a smooth curve germ 
(since $\bigl(g'_1(0),\dots, g'_{n-1}(0)\bigr)\neq (0,\dots,0)$)
which is simply tangent to  the branch locus $(y_n=0)$
(since $g_n(t)=f^2_n(t)$ vanishes with multiplicity 2 at $t=0$).

Let us now consider a complex analytic deformation
$$
\bigl(G_1(t,s),\dots, G_n(t,s)\bigr)\qtq{of}
\bigl(g_1(t),\dots, g_n(t)\bigr),
$$
where $s$ varies in a polydisc  $D^r$,
the $G_i$ are analytic  and
 $G_i(t,0)=g_i(t)$. When can we lift this local deformation of
$\pi\circ f$  to a deformation of $f$?
From (\ref{lifting.say}.2) we see that our only choice is
to take
$F_i(t,s)=G_i(t,s)$ for $i<n$ and
$F_n(t,s)=\sqrt{G_n(t,s)}$. That is,
the lifting is possible iff $ G_n(t,s)$
is a square. 

(There are two possible choices of the square root,
but only one of these
will agree with $f_n(t,0)$ when $s=0$.  Thus, if  a lifting exists,
it is unique.)
\medskip

{\it Lemma} \ref{lifting.say}.3. There is a hypersurface
$H\subset D^r$ such that 
$\bigl(G_1(t,s),\dots, G_n(t,s)\bigr)$ lifts to a deformation
$\bigl(F_1(t,s),\dots, F_n(t,s)\bigr)$ iff $s\in H$.
\medskip

Proof. By assumption, $G_n(t,s)$ contains $t^2$ with nonzero
coefficient. By the Weierstrass preparation theorem
(cf.\ \cite[p.8]{gri-har}) we can write
$$
G_n(t,s)= U(t,s)\cdot \bigl( t^2+b(s)t+c(s)\bigr)
$$
where $U(0,0)\neq 0$. Thus  $ G_n(t,s)$
is a square iff $b(s)^2-4c(s)=0$. Thus $H\subset D^r$
is defined by the equation $b(s)^2-4c(s)=0$.\qed
\medskip

{\it Corollary} \ref{lifting.say}.4. Composing with $\pi$
gives
$$
\pi_*: \mor(C,X)\to \mor(C, \p^n)\qtq{and}
\pi_*: \curves(X)\to \curves(\p^n).
$$
If $C$ satisfies the assumptions (\ref{lifting.say}.1) then
both of these are
 local embeddings near $[C]$ and 
their image has codimension $\leq (R\cdot C)$.
\medskip

Proof. If $\pi:X\to \p^n$ is a
local analytic isomorphism between neighborhoods $(x\in U_x)$ 
and $(\pi(x)\in V_x)$ then everything automatically lifts from
$V_x$ to $U_x$. Thus the only problem is
at the points $x\in R\cap C$. For every such point,
it is one condition to lift by (\ref{lifting.say}.3).
Thus there is a global lifting over the
intersection of the hypersurfaces $\{H_x: x\in R\cap C\}$. \qed
\end{say}

We are now ready to prove the main dimension estimates
for the spaces of curves. 
The formula is very much in the spirit of Riemann's original version of the
Riemann-Roch theorem: we estimate
$\dim \mor(C,X)$ from below in terms of intersection numbers
involving   Chern classes.
However, it is not known how to define appropriate analogs of
the higher cohomology groups that would make
the inequality into an equality,
and thus establish a better parallel with the
Riemann-Roch theorem.

\begin{thm}\label{dim.est.gen.thm}
 Let $X$ be a smooth quasi projective variety, $C$  a
proper reduced curve and $f:C\to X$ a morphism.
Then
 $$
\begin{array}{rcl}
\dim_{[f]} \mor(C,X)&\geq &
c_1(X)\cdot f_*[C] +\dim X\cdot \chi(\o_C)\\
& = & \deg_C f^*T_X +\dim X\cdot \chi(\o_C)\\
& = & H^0\bigl(C, f^*T_X\bigr)- H^1\bigl(C, f^*T_X\bigr).
\end{array}
$$
\end{thm}

Proof. Choose a general $\pi:X\to \p^n$
such that the assumptions (\ref{lifting.say}.1) hold
and the degree of $\pi\circ f$ is high enough.
By (\ref{dc.exmp.3}),
$$
\dim_{[\pi\circ f]} \mor(C,\p^n)\geq 
c_1(\p^n)\cdot (\pi\circ f)_*[C] +\dim X\cdot \chi(\o_C)
$$
and by (\ref{lifting.say}.4) this gives that
$$
\dim_{[f]} \mor(C,X)\geq 
\pi^*c_1(\p^n)\cdot f_*[C] +\dim X\cdot \chi(\o_C)-
\bigl(R\cdot f_*[C]\bigr).
$$
By the Hurwitz formula
$c_1(X)=\pi^*c_1(\p^n)-R$, giving the  inequality
in (\ref{dim.est.gen.thm}).

The first equality in the statement is clear since
$c_1(X)=c_1\bigl(\det T_X\bigr)$ and 
$$
\deg_C f^*T_X=\deg_C  f^*\det T_X=c_1(X)\cdot f_*[C]
$$
by the projection formula.
The last equality is just Riemann-Roch for curves.
\qed
\medskip

The same argument also gives
the dimension estimate for the space of curves:

\begin{thm}\label{dim.chow/hilb.thm}
 Let $X$ be a smooth projective variety and 
 $C\subset X$  a smooth curve. Then
$$
\dim_{[C]} \curves(X)\geq  c_1(X)\cdot [C] +
(\dim X -3)\cdot \bigl(1-g(C)\bigr),
$$
where $\curves(X)$ denotes either the
Chow variety or the Hilbert scheme of curves on $X$. \qed
\end{thm}

It is not hard to modify our arguments to see that (\ref{dim.chow/hilb.thm})
also holds if  $C$ is a reduced curve with locally smoothable
singularities, see \cite[II.1.14]{rc-book}.

\begin{say}[{\bf A philosophical claim and a challenge}]

The philosophical claim is that the estimates
(\ref{dim.est.gen.thm}) and (\ref{dim.chow/hilb.thm}) are optimal.
We will see that in many cases indeed equality holds,
thus in this weak sense they are  optimal.
More substantively, if we take a general almost complex perturbation
of the complex structure of $X$, then
(\ref{dim.est.gen.thm}) and (\ref{dim.chow/hilb.thm}) should
become  equalities at every ``interesting'' point of
 the space of pseudo-holomorphic curves.
(The formula frequently miscounts the dimension for curves
$f:C\to X$ for which $C\to f(C)$ is a multiple cover.
In some sense, we can ignore these if we study curves
on $X$. However, if one looks at families of curves on $X$,
such multiple covers naturally arise as limits of embedded curves.
For definitions and more details, see Section 5.)
One may thus claim that (\ref{dim.est.gen.thm}) and (\ref{dim.chow/hilb.thm})
correctly compute the dimension of the spaces of curves that
persist under small  almost complex perturbations of $X$,
but they fail to take into account the curves that
exist only ``accidentally.''

For instance, if $X^3_e\subset \p^4$ is a smooth hypersurface of degree
$e$, then we get that
$$
\dim_{C} \curves(X^3_e)\geq (5-e)\deg C,
$$
and the philosophical claim is that equality should hold.
In particular, for $e=\deg X^3_e\geq 6$ this
predicts that there are no curves at all on $X$!
Equivalently, every curve on a smooth, 3-dimensional 
hypersurface of degree $\geq 6$ is ``accidental.''

On the other hand, 
from the point of view of almost complex manifolds,
being an algebraic variety is an ``accident,'' and
on an algebraic variety there are
many curves. 
One of the biggest  challenges
of the theory of curves on varieties 
 is to explain  how to correct the formulas
(\ref{dim.est.gen.thm}) and (\ref{dim.chow/hilb.thm}) for algebraic varieties.
After all, our aim is to study algebraic varieties, not
almost complex manifolds.

\end{say}

\begin{aside} In general, from (\ref{dc.exmp.4}.1), 
we get the more precise formula
$$
\dim_{[f]} \mor_L(C,X)\geq 
c_1(X)\cdot f_*[C] +\dim X\cdot \chi(\o_C) -g(C)
+ \dim X\cdot h^1(C,L)- h^1(C,L^e).
$$
For instance, if $C\subset X^n\subset \p^{n+1}$ is embedded
by (a subsystem of) the canonical system, then
$$
\dim_{[f]} \mor_{K_C}(C,X)\geq 
c_1(X)\cdot \phi_*[C] +\dim X\cdot \chi(\o_C)-g(C)
+ \dim X.
$$
Thus if a hypersurface
contains such a curve, it contains a family of such curves
whose dimension is  at least $(\dim X)$--larger than one would expect. 
\end{aside}

\section{Free curves and uniruled varieties}{\ }

As in the Riemann-Roch theorem, the next step is to ask
when the  inequality in (\ref{dim.est.gen.thm}) is an equality.
The following result describes  essentially the only known general case
 when this holds and 
 the local structure of $ \mor(C,X)$ is fully understood.

\begin{thm} \label{dim.est.gen.cor}
 Let $X$ be a smooth quasi projective variety, $C$  a
proper reduced curve and $f:C\to X$ a morphism.
If $H^1\bigl(C, f^*T_X\bigr)=0$ then
$ \mor(C,X)$ is smooth at $[f]$ of dimension
$\deg_C f^*T_X +\dim X\cdot \chi(\o_C)$.
\end{thm}

Proof.  Assume first that $H^1(C,\o_C)=0$.
If $H^1\bigl(C, f^*T_X\bigr)=0$ then, by (\ref{dim.est.gen.thm}), 
$\dim_{[f]} \mor(C,X)\geq  h^0\bigl(C, f^*T_X\bigr)$.
On the other hand, by
(\ref{tangent.space.thm}), the tangent space of $\mor(C,X)$ at
$[f]$ is $H^0\bigl(C, f^*T_X\bigr)$.
The dimension of the tangent space is always at least
the dimension and equality holds only at smooth points. 

The case when $H^1(C,\o_C)\neq 0$ is harder since we have not
proved that in (\ref{tangent.space.thm}.2) the map
on the right is surjective. For a complete proof, see \cite[Sec.I.2]{rc-book}.
\qed
\medskip

The above   is a very useful result if there are many curves
on $X$ for which the condition $H^1\bigl(C, f^*T_X\bigr)=0$ holds.
Our next aim  is to get a feeling how frequently this happens.

\begin{exrc} Let $A$ be an Abelian variety.
Then $T_A$ is trivial, hence
$h^1\bigl(C, f^*T_A\bigr)=\dim A\cdot h^1(C,\o_C)$
for every $C$. Since an Abelian variety does not contain
any rational curves, we see that
$H^1\bigl(C, f^*T_A\bigr)$ is never zero.

(More precisely, if $\phi:\p^1\to A$ is a constant map,
then  $H^1\bigl(\p^1, f^*T_A\bigr)=0$. Here
$\mor(\p^1,A)$ consist only of constant maps,
so $\mor(\p^1,A)\cong A$ which is smooth of dimension $\dim A$.
While it may seem silly to think about such cases,
it is  useful to study maps from reducible curves
to a variety that may be constant on some of the irreducible components.
In using induction, we frequently need to make sure
that our formulas work for  constant maps as well.)
\end{exrc}

\begin{exrc} Let $C\subset \p^n$ be a smooth complete intersection 
of hypersurfaces of degrees  $d_1,\dots, d_{n-1}$.  Check that
$H^1\bigl(C, f^*T_{\p^n}\bigr)=0$ if and only if either $g(C)\leq 1$
or  $C\subset \p^n$ is canonically embedded.

Up to permuting the $d_i$, the first happens only in the cases
$$
(1,1,\dots,1),(2,1,\dots,1),(3,1,\dots,1),(2,2,1,\dots,1).
$$
The second possibility holds only in the cases $(4), (3,2), (2,2,2)$.
\end{exrc}

\begin{exrc} Let $f:C\to \p^n$ be a smooth curve
such that $\deg C\geq 2g(C)-1$. Then 
$H^1\bigl(C, f^*T_{\p^n}\bigr)=0$.
\end{exrc}

\begin{exmp} Let $X$ be a projective homogeneous space
(for instance, $\p^n$, a smooth quadric, a Grassmannian, \dots).
Then $T_X$ is generated by global sections, hence
$f^*T_X$ is also generated by global sections
for every $f:C\to X$. If $C\cong \p^1$, this implies that
$H^1\bigl(\p^1, f^*T_X\bigr)=0$.

Thus $\mor(\p^1, X)$ is everywhere smooth and of the expected dimension.
\end{exmp}

This might  be the only case where
$\mor(\p^1, X)$ is everywhere nice:

\begin{conj} \cite{cam-pet-MR1087244}\label{homog.space.char.conj}
Let $X$ be a smooth projective variety.
The following are equivalent:
\begin{enumerate}
\item  $H^1\bigl(\p^1, f^*T_X\bigr)=0$ for every $f:\p^1\to X$.
\item There is a morphism $p:X\to Y$ such that
\begin{enumerate}
\item $X\to Y$ is a locally trivial fiber bundle whose
fibers are projective homogeneous spaces under a linear algebraic group, and
\item every map $\p^1\to Y$ is constant.
\end{enumerate}
\end{enumerate}
\end{conj}

(Here we allow the uninteresting special case
when every map $\p^1\to X$ is constant
and we  take $Y=X$.)

We have basically established above that (\ref{homog.space.char.conj}.2)
implies (\ref{homog.space.char.conj}.1). 
The converse is known in low dimension; essentially as a consequence
of much stronger classification results.
The case of complete intersections is proved in \cite{pandharipande-2004}.

Note also that, at first sight, 
the blow up $B_p\p^2$ of $\p^2$ at a point seems like a counter example.
Indeed, $\aut\bigl(B_p\p^2\bigr)$ is transitive
away from the exceptional curve $E\subset B_p\p^2$, hence
$H^1\bigl(\p^1, f^*T_{B_p\p^2}\bigr)=0$, unless $f(\p^1)\subset E$.
The normal bundle of $E$ itself is $\o_E(-1)$, thus
the tangent bundle restricted to $E$ is
$\o_E(2)+\o_E(-1)$ and its
$H^1$ is zero.

However, if $f:\p^1\to E$ is a degree 2 map, then
$f^*T_{B_p\p^2}\cong \o_{\p^1}(4)+\o_{\p^1}(-2)$,
and its $H^1$ is nonzero.

\begin{defn} A morphism $f:\p^1\to X$ is called {\it free}
if $f^*T_X$ is generated by global sections.
Equivalently, if $H^1\bigl(\p^1, f^*T_X(-1)\bigr)=0$.
Informally, these are the rational curves that
can be deformed in every possible direction in $X$.

Note also that if  $f^*T_X$ is generated by global sections
over a nonempty open set then it is generated by global sections
everywhere.

A morphism $f:\p^1\to X$ is called {\it very free}
if $H^1\bigl(\p^1, f^*T_X(-2)\bigr)=0$.
Since every vector bundle on $\p^1$ is a direct sum of line bundles,
this is equivalent to saying 
 that 
   $$
   f^*T_X=\tsum_{i=1}^{\dim X} \o_{\p^1}(a_i) \qtq{with}  a_i\geq 1.
   $$
Informally, these are the rational curves that
can be deformed in every possible direction in $X$
while keeping a given point fixed.

\end{defn}

The following theorem says that, in some sense,
every rational curve on a variety is either free or special.

\begin{thm} \label{vg.curve,free.thm}
Let $X$ be a smooth projective variety. 
Then there are countably many proper closed subvarieties
$V_i\subset X$ such that every
$f:\p^1\to X$
 whose image is not contained in any of the $V_i$ is free.
\end{thm}

Proof. There are countably many irreducible components
of $\mor(\p^1,X)$. Let us see how each of them leads
to finitely many of the $V_i$s.

To be slightly more general, let
$U\subset \mor(\p^1,X)$ be any irreducible subset
that we are interested in.
Consider the universal morphism
$\univ_U:U\times \p^1\to X$.
If $\univ_U$ is not dominant,  the closure of
$\univ_U\bigl(U\times \p^1\bigr)$ will be one of the $V_i$s. 

Thus assume that $\univ_U$ is dominant.
Then there is a dense open subset 
   $W\subset \p^1\times U$ such that $\univ_U:W\to X$ is smooth 
and its derivative
   $d(\univ_U): T_{\p^1\times U}\to \univ_U^*T_X$ is surjective over $W$.

   Let $\pi_2:\p^1\times U\to U$ be the second projection
   and assume that $u \in \pi_2(W)$. Then 
   $$
   \o_{\p^1}(2) \oplus \o_{\p^1}^{\dim U}\cong 
   T_{\p^1\times U}|_{\p^1\times\{u\}} 
   \to f_u^*T_X 
   $$
   is surjective on the open set
$W\cap \p^1\times\{u\}$. Therefore $f_u$ is free.

This takes care of the maps corresponding to points in
$\pi_2(W)$. The complement $U\setminus \pi_2(W)$
is a union of lower dimensional subvarieties.
We repeat the above argument for each, and so on.
Eventually, for any irreducible component $U_i\subset \mor(\p^1,X)$
 we get finitely many    proper closed subvarieties
$V_{ij}\subset X$ such that for every $[f]\in U_i$, either
$f(\p^1)$
 is  contained in some $V_{ij}$  or $f$ is free.

Since $ \mor(\p^1,X)$ has countably many  irreducible 
components, we may need to exclude  countably many 
subvarieties $V_{ij}$.
(It is very poorly understood when one actually needs
  countably many exceptions. Such examples are
given by $\p^2$ blown up at $m\geq 9$ general points
and by general K3 surfaces.)
\qed
\medskip

This shows that the 3 conditions in the following definition are
equivalent.

\begin{defn} A smooth projective variety $X$ over $\c$ is
{\it uniruled} iff the following equivalent  conditions hold.
\begin{enumerate}
\item There is a dense open set $X^0\subset X$ such that for 
   every $x\in X^0$, there is a rational curve through $x$.
\item There is a variety $Y$ and a
dominant and generically finite map $Y\times \p^1\map X$.
\item There is a  free morphism $f:\p^1 \to X$.
\end{enumerate}
\end{defn}

Thus, if $X$ is not uniruled, then all rational curves
on $X$ lie in the union of at most countably many
subvarieties $V_i\subsetneq X$. We think of this as having
only ``few'' rational curves on $X$.

\section{Very free curves and rationally connected varieties}{\ }

This section is mostly a summary of the basic results.
For a general overview, see \cite{k-bull}.
An introduction with proofs is given in \cite{ar-ko}, while
those wishing to go through all the details should consult the
original papers or \cite{rc-book}.

Being uniruled is not a good structural property,
since $Y\times \p^1$ is uniruled for any $Y$.
For over a century it has been an open problem to define a good subclass
of uniruled varieties that does not contain any such 
``mongrel''  examples and generalizes to higher dimensions
the following result about surfaces

\begin{say}[Basic trichotomy of surface theory] \label{trich.surf.say}
Let $S$ be a smooth
projective surface. Then exactly one of the following holds:
\begin{enumerate}
\item  $S$ is not uniruled (and hence has at most 
countably many rational curves),
\item $S$ is uniruled but maps to a non-uniruled 
(equivalently, of genus $\geq 1$) curve 
(and hence all rational curves on $S$ are in the fibers of this map), or
\item $S$ is birational to $\p^2$ (and hence has many rational curves).
\end{enumerate}
\end{say}

By now it is clear that
the key property is to require rational curves not just through
any point, but through any pair of points. One can then imagine several
variants of this concept. The next result shows that
all of these are equivalent.
For the proofs, see \cite{kmm2, rc-book, ar-ko, sankaran-2007}.

\begin{thm} \label{RCV.thm} Let $X$ be a smooth projective variety over
$\c$. The
following are equivalent.

\begin{enumerate}
\item There is a dense open set $X^0\subset X$ such that for 
   every $x_1,x_2\in X^0$, there is a rational curve through $x_1$ and $x_2$.
\item For every $x_1,x_2\in X$, there is a rational curve 
   through $x_1$ and $x_2$.
\item For every integer $m>0$, and every $x_1,\dots ,x_m\in X$, 
   there is a rational curve through $x_1,\dots ,x_m$.
\item For every 0-dimensional subscheme $Z\subset \p^1$,
every morphism $f_Z:Z\to X$ can be extended to a morphism
$f:\p^1\to X$.
\item There is a dense open set $X^0\subset X$ such that, for 
   every $x_1,x_2\in X^0$, there is a connected curve with rational components
   through $x_1$ and $x_2$.
\item There is a very free morphism $f:\p^1 \to X$.
\item Let $C$ be a smooth curve, 
$Z\subset C$ be a 0-dimensional subscheme,  $r>0$ an integer
and $W\subset X$ a subscheme of codimension $\geq 2$. Then 
every morphism $f_Z:Z\to X$ can be extended to a  morphism
$f:C\to X$ such that
\begin{enumerate}
\item $H^1\bigl(C, f^*T_X(-r)\bigr)=0$, 
\item $f\bigl(C\setminus Z\bigr)$ is disjoint from $W$, 
\item $f$ is an embedding on $C\setminus Z$ if $\dim X\geq 3$,
\item $f$ is an embedding on $C$ if $\dim X\geq 3$ and
$f_Z$ is an embedding.
\end{enumerate}
\end{enumerate}
\end{thm}

\begin{defn} \label{RCV.defn} Let $X$ be a smooth projective variety over
$\c$. We say that $X$ is {\it rationally connected} if it satisfies 
the equivalent conditions of (\ref{RCV.thm}).
\end{defn}

There are two additional (partially conjectural) characterizations that
are of interest.

\begin{conj} Let $X$ be a smooth projective rationally connected
variety. Let $C$ be a smooth curve, $D\subset C$  a Euclidean open set and
$\phi_D:D\to X$ a holomorphic map.
Then
there is a sequence of algebraic maps $f_r:C\to X$
such that the $f_r|_D$ converge to $\phi_D$
in the compact-open topology.

Moreover, if $Z\subset D$ is a finite set then we can
also assume that $f_r|_Z=\phi_D|_Z$.
\end{conj}

\begin{say}[Loop spaces of rationally connected varieties]
\cite{lem-sza} proves that a variety $X$ is rationally
connected iff its loop space $\Omega X$   is rationally
connected. By a loop space we mean the space of all
 continuous/differentiable maps form the circle $S^1$ to $X$.
(The loop spaces have  a natural complex Banach  manifold
structure, see \cite{lem-sza} for detail.)
Aside from technicalities, the key result is the
following, which fits very nicely in the sequence of
characterizations in (\ref{RCV.thm}).
\begin{enumerate}\setcounter{enumi}{7}
\item Fix $m>0$  continuous/differentiable maps
$\phi_i:S^1\to X$ and distinct points $p_i\in \p^1$.
Then
there is a sequence of   continuous/differentiable maps
$\Phi_r:S^1\times \c\p^1$ that are algebraic on each 
$\{s\}\times \c\p^1$ such that, as $r\to\infty$, the 
$\Phi_r|_{S^1\times \{p_i\}}$ converge to $\phi_i$ for every $i$.
\end{enumerate}

Conjecturally, one can even achieve that
$\Phi_r|_{S^1\times \{p_i\}}=\phi_i$ for every $i$,
but this is proved only for general $\phi_i$.
\end{say}

Being rationally connected is stable under various operations:

\begin{thm} \cite{kmm2} \label{RC.stable.props}
Let $X$ be a smooth projective rationally connected variety.
Then:
\begin{enumerate}
\item Every  smooth projective  variety that is birational to $X$
is also rationally connected.
\item Every  smooth projective  variety that is the image of $X$
by a rational map 
is also rationally connected.
\item Every  smooth projective  variety that is a deformation of $X$ 
is also rationally connected.
\end{enumerate}
\end{thm}

It is also useful to know that many varieties are rationally connected:

\begin{thm} \cite{nadel, kmm1, kmm3, camp}
 Let $X$ be a smooth, projective  Fano variety,
that is, with $-K_X$ ample. Then $X$ is rationally connected.
\end{thm}

\begin{thm}\cite{ghs} Let $X$ be a smooth, projective variety
and $f:X\to Y$ a morphism. Assume that $Y$ is rationally connected
and so  are the smooth fibers of $f$. 
Then $X$ is also rationally connected.
\end{thm}

We also have the basic trichotomy of algebraic varieties,
which is a close analog of (\ref{trich.surf.say}).
The map in (\ref{trich.gen.say}.2) is called the
{\it MRC-fibration} or {\it maximal rationally connected fibration} of $X$.

\begin{thm}\cite{kmm2}  \label{trich.gen.say} Let $X$ be a smooth
projective variety. Then exactly one of the following holds:
\begin{enumerate}
\item  $X$ is not uniruled (and hence has  ``few'' rational curves),
\item there is a map $\pi:X\map Y$ onto  a non-uniruled variety
with $0<\dim Y <\dim X$
whose general fibers are rationally connected 
(and hence most rational curves on $X$ are in the fibers of $\pi$), or
\item $X$ is rationally connected (and hence has many rational curves).
\end{enumerate}
\end{thm}

\section{Connections with symplectic geometry}{\ }

For general introductions to the topics in this section, see
\cite{ful-pan, mcd-sal1, mcd-sal2} and
for  more technical details consult
\cite{beh-fan, li-ti1, li-ti2}.

\begin{say}[Symplectic structure of varieties]

Any smooth projective variety admits a symplectic structure. This can be
constructed as follows. 
  On $\c^{n+1}$
consider the Fubini--Study 2-form
$$
\omega':=
\frac{\sqrt{-1}}{2\pi}
\left[
\frac{\sum dz_i\wedge d\bar z_i}
{\sum |z_i|^2}
-
\frac{(\sum \bar z_i dz_i)\wedge (\sum  z_i d\bar z_i)}
{\left(\sum |z_i|^2\right)^2}
\right].
$$
It is closed, non-degenerate on $\c^{n+1}\setminus\{0\}$ and invariant under
scalar multiplication. Thus $\omega'$ descends to a symplectic 2-form
$\omega$ on $\c\p^n=(\c^{n+1}\setminus\{0\})/\c^*$. 
This construction depends on the choice of a basis in $\c^{n+1}$ and
one can see that $\omega$ is invariant under the unitary group
$U(n+1)$ but not under $\aut(\c\p^n)=PGL(n+1)$.
Thus it is better to think of $\c\p^n$ yielding
not just one  symplectic manifold  $\bigl(\c\p^n,\omega\bigr)$
but rather a whole family of  symplectic manifolds
parametrized by the connected space $PGL(n+1)/U(n+1)$.
 
Generalizing this, 
we say that two symplectic manifolds $(M,\omega_0)$ and $(M,\omega_1)$ are
{\it symplectic deformation equivalent} if there is a continuous family of 
symplectic manifolds $(M,\omega_t)$ starting with
$(M,\omega_0)$ and ending with $(M,\omega_1)$.

If $X\subset \c\p^n$ is any smooth variety, then the restriction
$\omega|_X$ makes $X(\c)$  into a symplectic manifold. 
(Note that this $\omega|_X$ has nothing to do with the
dualizing sheaf, commonly denoted by  $\omega_X$.)
The resulting symplectic manifold $(X(\c), \omega|_X)$ depends
on the choice of $\omega$, but the dependence is rather
clear.  
Thus to every smooth projective variety 
and the choice of a (very) ample cohomology class,
the above construction associates a
symplectic manifold $(X(\c), \omega|_X)$ which is unique up to symplectic
deformation equivalence.
\end{say}

A powerful way to relate properties of
the symplectic manifold $(X(\c), \omega|_X)$
to the algebraic geometry of $X$ is through the 
enumerative properties 
of stable curves.

\begin{defn} Let $X$ be a variety. A genus $g$ stable
 stable  curve with $n$ marked points  over $X$
 is a triple $(C,P,f)$,
where 
\begin{enumerate}
\item $C$ is a proper connected curve  having only
nodes,  
\item $P=(p_1,\dots, p_n)\subset C$ is an ordered set of smooth
points of $C$, 
\item  $f:C\to X$ is a morphism, and 
\item $C$ has only finitely many automorphisms that 
fix $P$ and commute with $f$.
\end{enumerate}
As for any map of a curve to a variety,
$f_*[C]\in H_2(X(\c),\z)$ is well defined. 

For a given
$\beta\in H_2(X(\c),\z)$ and subvarieties
$Z_i\subset X$,  let
$$
{\mathcal M}_{g,n}(X,\beta, Z_1,\dots, Z_n)
$$ denote the set of all
genus $g$ stable
  curves with $n$ marked points such that $f_*[C]=\beta$
and $f(p_i)\in Z_i$ for $i=1,\dots,n$.
\end{defn}

\begin{say}[Gromov-Witten invariants]\label{GW.say}
It is not hard to define families of stable curves,
 interpret ${\mathcal M}_{g,n}(X,\beta, Z_1,\dots, Z_n)$ as a 
 moduli functor and define its coarse moduli space
$ M_{g,n}(X,\beta, Z_1,\dots, Z_n)$. 
We are mostly interested in the cases when
$ M_{g,n}(X,\beta, Z_1,\dots, Z_n)$ is a finite set
with reduced scheme structure. In these cases 
functorial  definitions add nothing interesting to the picture.

An easy  generalization of  the dimension estimates
in Section 2  shows that
$$
\begin{array}{l}
\dim M_{g,n}(X,\beta, Z_1,\dots, Z_n)\geq\\
\qquad\qquad \geq 
c_1(X)\cdot \beta +(\dim X-3)\cdot \chi(\o_C) -\sum_i (\codim_XZ_i-1),
\end{array}
\eqno{(\ref{GW.say}.1)}
$$
(There are two ways to think about this formula.
First, in the third step of (\ref{philos.say.1}) we can
use $M_g$. Then the last term in (\ref{GW.say}.1) 
 comes from the observation that
in a family of curves on $X$ it is $(\codim_XZ-1)$ conditions
 to intersect a subvariety
$Z\subset X$.
Alternatively, we can use $ M_{g,n}$  (which adds an extra $n$
to the formula) and note that in a family of pointed 
curves on $X$ it is $\codim_XZ$ conditions
 for the marked  point to lie on 
$Z$. Since there are $n$ of the $Z_i$, this
cancels out the extra $n$.)

If the right hand side of (\ref{GW.say}.1) is 0, we expect
that there are only finitely many maps
in  ${\mathcal M}_{g,n}(X,\beta, Z_1,\dots, Z_n)$
and their number, called a {\it Gromov-Witten invariant} of $X$,
 has  enumerative significance.

In fact, one can define a Gromov-Witten invariant
$$
\Phi_{g,n}(X,\beta, Z_1,\dots, Z_n)
\eqno{(\ref{GW.say}.2)}
$$
whenever the right hand side of (\ref{GW.say}.1) is 0,
even if $ M_{g,n}(X,\beta, Z_1,\dots, Z_n)$ is positive dimensional.
One  advantage of  Gromov-Witten invariants is that they 
depend on the 
symplectic structure only. That is, although 
$M_{g,n}(X,\beta, Z_1,\dots, Z_n)$
visibly depends on the algebraic variety $X$ and the $Z_i$, the 
value of the Gromov-Witten invariants
depends only  on the symplectic deformation class of $(M_X,\omega)$, 
the numbers $g,n$
and the homology classes $\beta, [Z_i]$:
$$
\Phi_{g,n}(X,\beta, Z_1,\dots, Z_n)=
\Phi_{g,n}\bigl((M_X,\omega), \beta, [Z_1],\dots, [Z_n]\bigr).
\eqno{(\ref{GW.say}.3)}
$$
In general,   Gromov-Witten invariants are rational numbers
which can be negative or zero.  
(A curve $C$ with automorphisms may count as $1/|\aut(C)|$
and positive dimensional components may have negative
contribution; see (\ref{rahul.exmp}) for such an example.)
There are, however,   a few cases when the obvious algebraic count
gives the  Gromov-Witten invariant.
\medskip

{\it Example} \ref{GW.say}.4. Assume that 
${\mathcal M}_{g,r}(X,\beta, Z_1,\dots, Z_r)=\emptyset$.
Then,  for any homology classes $[Z_{r+1}],\dots, [Z_n]$,
the Gromov-Witten invariant 
$\Phi_{g,n}(M_X,\omega, \beta, [Z_1],\dots, [Z_n])$ is zero.

\medskip

{\it Example} \ref{GW.say}.5. Assume that for  every
$(C,f)\in  {\mathcal M}_{0,0}(X,\beta)$ we have that 
$C\cong \p^1$, $f$  is an immersion
and  $f^*T_X=\o(2)+\o(-1)^{\dim X-1}$. Then
the Gromov-Witten invariant equals the number of such maps:
$$
\Phi_{0,0}(X,\beta)=\# {\mathcal M}_{0,0}(X,\beta).
$$
\medskip

{\it Example} \ref{GW.say}.6. The case when  the $Z_i$
are ``sufficiently general,'' can be reduced to
(\ref{GW.say}.5). First we need to assume that,
for every $i$, $f(p_i)$ is a smooth point of $Z_i$
and $p_i= f^{-1}(Z_i)$  scheme theoretically
(that is, $Z_i$ and $f(C)$ meet only at $f(p_i)$ and transversally).

If this holds, let $B_ZX\to X$ denote the blow up of $\cup Z_i\subset X$,
 $\tilde f:C\to B_ZX$ the lifting of $f$
and $\beta_Z:=\tilde f_*[C]$.
($B_ZX$  can be a quite singular space. However, by our assumptions,
$\cup Z_i$ is smooth in a  neighborhood of $f(C)$,
thus the singularities of $B_ZX$ will not matter to us.)

Under these assumptions, $f\mapsto\tilde f$
establishes a bijection 
$$
{\mathcal M}_{g,n}(X,\beta, Z_1,\dots, Z_n)
\quad\leftrightarrow\quad
{\mathcal M}_{g,0}(B_ZX,\beta_Z).
$$
 That is, the blow up removed the subvarieties
$Z_i$ and the marked points $P\subset C$ from the picture.


Now assume further
 that, for  every
$(C,P, f)\in  {\mathcal M}_{0,n}(X,\beta, Z_1,\dots, Z_n)$, the curve 
$C$ is isomorphic to $\p^1$,  $f$  is an embedding
and  $\tilde f^*T_X=\o(2)+\o(-1)^{\dim X-1}$. Then
the Gromov-Witten invariant is just the number of such maps:
$$
\Phi_{0,n}(X,\beta, Z_1,\dots, Z_n)=
\#{\mathcal M}_{0,n}(X,\beta, Z_1,\dots, Z_n).
$$
\medskip

{\it Example} \ref{GW.say}.7. More generally, assume that the
expected and the true dimensions of ${\mathcal M}_{g,n}(X,\beta, Z_1,\dots, Z_n)$
are both  zero and the corresponding maps are birational to their image. Then
 $$
\Phi_{g,n}(X,\beta, [Z_1],\dots, [Z_n])\geq 
\#{\mathcal M}_{g,n}(X,\beta, Z_1,\dots, Z_n).
$$


\end{say}

We use these examples to  show that being uniruled
is a property of the underlying symplectic variety.

\begin{thm}\cite[4.2.10]{k-low}, \cite{ruan-1}\label{uniruled.symp.thm}
 Let $X_1, X_2$ be smooth projective varieties and
 $\bigl(M_i, \omega_i\bigr)$ the corresponding symplectic manifolds.
Assume that $\bigl(M_1, \omega_1\bigr)$
is symplectic deformation equivalent to $\bigl(M_2, \omega_2\bigr)$.
Then $X_1$ is uniruled iff $X_2$ is.
\end{thm}

Proof. 
Assume that $X_1$ is uniruled. 
Fix a very general point $x_1\in X_1$
such that the corresponding point $x_2\in X_2$ is also very general. 
 Let $H$ be  a very ample divisor on $X_1$ and 
 $f:\p^1\to X_1$  a map such that $f(0)=x_1$ and
$\bigl(H\cdot f_*[\p^1]\bigr)$ is the smallest possible.
Set $\beta:=f_*[\p^1]$ and 
consider  ${\mathcal M}_{0,1}\bigl(X_1,\beta, x_1\bigr)$.
Since $\bigl(H\cdot f_*[\p^1]\bigr)$ 
is the smallest possible, every such curve is irreducible.
Since $x_1\in X_1$ is very general,
every such curve 
is free (\ref{vg.curve,free.thm}). Set
$r:=\dim {\mathcal M}_{0,1}\bigl(X_1,\beta, x_1\bigr)$, 
let $H_1,\dots, H_{2r}\subset X_1$ be general divisors linearly
equivalent to $H$ and consider
$$
{\mathcal M}_{0,r+1}
\bigl(X_1,\beta, x_1, H_1\cap H_2,\dots, H_{2r-1}\cap H_{2r}\bigr),
$$
that is, stable rational curves in the homology class $\beta$
that pass through $x_1$ and intersect each of the
$H_{2j-1}\cap H_{2j}$.

The assumptions of  (\ref{GW.say}.7)
are satisfied and therefore
$$
\begin{array}{l}
\Phi_{0,r+1}
\bigl(X_1,\beta, x_1, [H_1\cap H_2],\dots, [H_{2r-1}\cap H_{2r}]\bigr)\geq \\
\qquad\qquad\qquad \geq 
\# {\mathcal M}_{0,r+1}
\bigl(X_1,\beta, x_1, H_1\cap H_2,\dots, H_{2r-1}\cap H_{2r}\bigr)>0.
\end{array}
$$
(By  \cite[II.3.14.4]{rc-book}, all our maps are immersions,
thus we could have used the more elementary  (\ref{GW.say}.6).)

Since the Gromov-Witten numbers are symplectic invariants,
this implies that 
$$
\Phi_{0,r+1}
\bigl(X_2,\beta, x_2, [H_1\cap H_2],\dots, [H_{2r-1}\cap H_{2r}]\bigr)>0.
$$
Thus, by (\ref{GW.say}.4), 
${\mathcal M}_{0,1}\bigl(X_2,\beta, x_2\bigr)$
is not empty and so $X_2$ also contains a rational curve
through $x_2$.  
(Note that  ${\mathcal M}_{0,1}\bigl(X_2,\beta, x_2\bigr)$
may not
contain any irreducible curves, but all the irreducible components
of the curves in ${\mathcal M}_{0,1}\bigl(X_2,\beta, x_2\bigr)$
are rational, so we do get an irreducible rational curve
through $x_2$ after all.) Thus $X_2$ is also uniruled. 
\qed

\medskip

This was  only  a warm-up for the
main question, which would be a significant 
generalization of (\ref{RC.stable.props}.3).

\begin{conj}\cite[4.2.7]{k-low} \label{rc.symp.conj}
Let $X_1, X_2$ be smooth projective varieties and
 $\bigl(M_i, \omega_i\bigr)$ the corresponding symplectic manifolds.
Assume that $\bigl(M_1, \omega_1\bigr)$
is symplectic deformation equivalent to $\bigl(M_2, \omega_2\bigr)$.
Then $X_1$ is rationally connected  iff $X_2$ is.
\end{conj}

This is still wide open, even for 3-folds, despite 
significant partial results  in \cite{voisin}.

\begin{say} Instead of discussing the positive results, let us illustrate
the problem  by a simple example showing 
why the proof of (\ref{uniruled.symp.thm})
does not work directly in the  rationally connected case.

For $e\geq 0$ consider the minimal ruled surface
$\pi:\f_{2e}\to \p^1$. The class of a  fiber is denoted by $F$
and the negative section by $E_{2e}$. Then
$H_2(\f_{2e},\z)=\z[E_{2e}]+\z[F]$, 
$(E_{2e}^2)=-2e, (E_{2e}\cdot F)=1$ and $(F^2)=0$.  

Note that $\f_{2e}$ is deformation equivalent to
$\f_0\cong \p^1\times \p^1$. Under this equivalence,
the fibers correspond to each other but
$E_{2e}$  corresponds to $E_0-eF$.

Given 2 general points $p,q\in \f_{2e}$, the smallest degree
curve (in any projective embedding)
connecting $p$ and $q$ is a reducible curve,
consisting of the two fibers  $F_p$ (resp.\ $F_q$) through $p$ (resp.\ $q$)
 and of $E_{2e}$.
Its homology class is  $\bigl[E_{2e}+2F\bigr]$
and $E_{2e}+F_p+F_q$ is the unique curve passing through
$p,q$ whose homology class is  $\bigl[E_{2e}+2F\bigr]$. 

We may be tempted to believe that 
having such a curve
 is a symplectic property.
However, this fails,  even for smooth algebraic deformations. 
Indeed, 
by the above remarks, 
the homology class $\bigl[E_{2e}+2F\bigr]$ becomes
$\bigl[E_0+(2-e)F\bigr]$ on $\f_0$, and 
there is no effective curve
on $\f_0$ in the homology class $\bigl[E_0+(2-e)F\bigr]$ if $e>2$.

We can try next to work with irreducible curves.
The smallest homology class that contains an irreducible rational curve 
through $p,q$ is $\bigl[E_{2e}+2eF\bigr]$.
The linear system $\bigl|E_{2e}+2eF\bigr|$ has dimension
$2e+1$ and its general member is a smooth rational curve. 
Those curves that pass through $p,q$ form a  linear subsystem of 
dimension $2e-1$. At first sight it seems that we can repeat
the arguments of (\ref{uniruled.symp.thm}).

There is, however, a hitch. We need to consider not the
space of curves in $\f_{2e}$ but the space of maps of curves
to  $\f_{2e}$. For irreducible curves, these two spaces are essentially
the same, but in general problems arise with multiple covers.
In our case, we can have maps whose set-theoretic image is 
 $F_p+F_q+E_{2e}$, but give a degree $2e-1$ cover over $F_p$.
These have a moduli space of dimension $2(2e-1)-2=4e-4$.
Thus, for $e\geq 2$, this has greater dimension
than the ``main component'' which is birational
to the linear system  $\bigl|E_{2e}+2eF\bigr|(-p-q)$.

As we see in (\ref{rahul.exmp}), 
larger dimensional components may give a negative contribution
to a Gromov-Witten invariant which may cancel the positive
contribution given by the irreducible curves.

Actually, in this case, we end up with a correct argument
if we follow the method of (\ref{uniruled.symp.thm}).

By fixing $2e-1$ other general points
$r_1,\dots, r_{2e-1}$, we see that there is a unique
curve in  $\bigl|E_{2e}+2eF\bigr|$ passing through all the points
$p,q,r_1,\dots, r_{2e-1}$. This curve is smooth and irreducible
and gives $1$ for the value of the Gromov-Witten invariant.
The larger dimensional components do not contribute anything.
This, however, seems more luck than a general principle.
\end{say}

The following example, based on a suggestion of R.~Pandharipande, 
illustrates that negative contributions 
from too large components can cancel out a nice curve, even
for algebraic deformations.

\begin{exmp}\label{rahul.exmp}
In $\p^3$ consider the family of quadrics
$$
Q_t:=\bigl(x_0^2-x_1x_2-t^2x_3^2=0\bigr).
$$
For $t\neq 0$ we get a smooth quadric, isomorphic to
$\f_0\cong \p^1\times \p^1$. 
Let $|L_t|$ denote one of the two families of lines.
For $t=0$ we get a singular quadric. We can resolve the singularity
by blowing up $(x_0-tx_3=x_1=0)$. We get a family of smooth
surfaces. For $t\neq 0$ we still have $\f_0$, but for
$t=0$ we get $\f_2$, the blow-up of $Q_0$ at the origin.
Let $E_2\subset \f_2$ be the exceptional curve and
$|L_0|$ the birational transforms of the family of lines on $Q_0$.
Note that as $t\to 0$, the limit of the family of lines
 $|L_t|$ is the family $E_2+|L_0|$ of reducible curves.

Let $|H|$ be the  pull-back of the family of hyperplane sections
of $Q_0$ to $\f_2$. Its 
singular members  form the family
  $E_0+2|L_0|$. 

Over the
pair of lines $(st=0)\subset \c^2_{s,t}$ 
consider a family of curves and smooth surfaces 
as follows.
\begin{enumerate}
\item Over the $t$-axis, we have the family $Q_t$ degenerating to $\f_2$
and curves $2L_t$  degenerating to  $2E_2+2L_0$.
\item Over the $s$-axis, we have the trivial family of  $\f_2$
with  curves  $E_2+H_s$  degenerating to  $2E_2+2L_0$.
\end{enumerate}

Set $X:=\f_2\times \p^1$,
let $\beta$ be the class of $\{(\mbox{point})\}\times \p^1$
 and set $Z_i:=(E_2+H^i)\times \{p_i\}$ where
$p_i\in \p^1$ are two distinct points and
 $H^i\in|H|$ are two smooth hyperplane sections intersecting
at 2  distinct points $q_1, q_2\in \f_2$.

Consider ${\mathcal M}_{0,2}(X, \beta, Z_1, Z_2)$. 
It consist of two isolated points corresponding to the curves 
$\{q_i\}\times \p^1\into \f_2\times \p^1$
 (each contributing 1 to the Gromov-Witten invariant)
and a 1-dimensional component
of curves of the form $\{(\mbox{point})\}\times \p^1$ that are contained
 in $E_2\times \p^1$. The expected dimension is  0.

If we move over to $\f_0\times \p^1$,
the $Z_i$ can be represented by curves
of the form $L_i+L'_i$ where
 the lines $L_i,L'_i$ are in the same family
of lines. Thus we can choose $Z_1$and $Z_2$  to be disjoint,
showing that  ${\mathcal M}_{0,2}(X, \beta, Z_1, Z_2)=\emptyset$.

\end{exmp}

 \begin{ack} I thank 
A.~Cannas~da~Silva for the invitation, the 
great organization and hospitality and 
E.~Ionel, S.~Katz, J.~Li, R.~Pandharipande, J.~Starr, C.~Voisin 
and the referees for useful comments, corrections and references.
Partial financial support  was provided by  the NSF under grant number 
DMS-0758275.
\end{ack}

\def\cprime{$'$} \def\cprime{$'$} \def\cprime{$'$} \def\cprime{$'$}
\providecommand{\bysame}{\leavevmode\hbox to3em{\hrulefill}\thinspace}
\providecommand{\MR}{\relax\ifhmode\unskip\space\fi MR }
\providecommand{\MRhref}[2]{%
  \href{http://www.ams.org/mathscinet-getitem?mr=#1}{#2}
}
\providecommand{\href}[2]{#2}


\noindent Princeton University, Princeton NJ 08544

\begin{verbatim}kollar@math.princeton.edu
\end{verbatim}

\end{document}